\documentclass[12pt]{amsart}

\usepackage{amsfonts}
\usepackage{amsthm}
\usepackage{amsmath}
\usepackage{amssymb}
\usepackage[cp1251]{inputenc}
\usepackage[english]{babel}
\newcommand {\bpr} {\noindent{ P r o o f.} }
\theoremstyle{plain}
\newtheorem{Th}{Theorem}
\newtheorem{Le}{Lemma}
\newtheorem{Pro}{Proposition}
\newtheorem{Def}{Definition}
\def\bs{~\hfill\rule{7pt}{7pt}}

\def\N{\mathbb N}

\begin{document}

\title{Almost periodic discrete sets}

\author{Favorov S., Kolbasina Ye.}

\address{Mathematical School, Kharkov National University, Swobody sq.4,
Kharkov, 61077 Ukraine}

 \email{Sergey.Ju.Favorov@univer.kharkov.ua,kvr${\_}$jenya@mail.ru}

\date{}
\begin{abstract}
 Using a special metric in the space of sequences, we give a geometric description of
almost periodic sets in the $k$-dimensional Euclidean space. We prove the completeness
of the space of almost periodic sets and some analogue of the Bochner criterion of
almost periodicity.
 Also, we show the connection between these sets and almost periodic measures.
\end{abstract}

\keywords{almost periodic measure, almost periodic set, almost periodic mapping}

\subjclass{Primary: 54E35; Secondary:52C23, 11K70}

\maketitle

 {\bf Introduction.} Almost periodic (with respect to the real shifts) sets were first
regarded by M.G.Krein and B.Ya.Levin \cite{krein lev} (see also the monograph
\cite{lev}) when studying the zero distribution of entire almost periodic functions.

\noindent{\bf Definition.}
A discrete set $\{a_n\}\subset \{ z: |Im z |< M <\infty\}$ is called almost periodic if for each
$\varepsilon >0$ there exists $ L < \infty $ such that every real interval of length $ L $ contains at
least one number $ \tau $ such that for some bijection $\rho
:\mathbb{Z}\rightarrow\mathbb{Z}$
 $$
 |(a_n+\tau)-a_{\rho(n)}|<\varepsilon\quad \textrm{ for all } \,n\in\mathbb{Z}.
  $$

\noindent M.G.Krein and B.Ya.Levin have shown that the zero set of an almost periodic
function from some special class (the class $[\Delta]$, see \cite{krein lev}) is an
almost periodic set. Moreover, any almost periodic set in a strip can be completely
described as the zero set of some analytic function with almost periodic modulus (see
\cite{rrf}).

In the present paper we study sets in $ \mathbb{R}^k$ which are almost periodic with
respect to  arbitrary shifts. Such sets arise in the models describing quasicristallic
structures (see, for example, \cite{lag}).

Using a special metric in the space of sequences, we give a geometric description of
almost periodic sets. We prove the completeness of the space of almost periodic sets
and some analogue of the Bochner criterion of almost periodicity. Following \cite{ron
tube}, we introduce a notion of almost periodic measure in $ \mathbb{R}^k$. We show
that a set is almost periodic if and only if the discrete measure with unit masses at
the points of this set is almost periodic.

 Throughout the work we denote $ i$-coordinate of a point $ x \in\mathbb{R}^k $ by $ x^i $, an open
 $k$-dimensional ball of radius $ R $ with center $ x\in\mathbb{R}^k $ by $ B(x,R)$, $k$-dimensional
 cube $ \{ y \in \mathbb{R}^k \mid x^i-L/2 \leq y^i < x^i+L/2, i=\overline{1,k}\}$ by $ Q(x,L) $. For
 any set $ A \subset\mathbb{R}^k $ for any $ \rho>0$ we write $ \; A_\rho=\displaystyle\bigcup_{x\in A}B(x,\rho)$
 and $ A_{-\rho} =\{ x \in \mathbb{R}^k \mid B(x,\rho) \subset A\}$. We denote the interior of $A$ by $Int A$,
 $k$-dimensional Lebesgue measure of $ A $ by $ m(A) $, the inner product of elements $x,y \in\mathbb{R}^k $
 by $\langle x,y\rangle$,  usual Euclidean norm by $|x|$. For a mapping  $g(x): \mathbb{R}^k \longrightarrow
 \mathbb{R}^n$ and $ \tau \in \mathbb{R}^k $
 we put $g^{\tau}(x)= g(x-\tau) $, for a measure $ \mu $ on $ \mathbb{R}^k $ we put
 $\mu^{\tau}(E)=\mu(E +\tau)$ for any Borel set $E\subset\mathbb{R}^k$.

\vspace{10pt}

{\bf 1. Almost periodic mappings.} Call to mind some definitions concerning almost
periodicity.

Let $ g $ be a continuous mapping from $ \mathbb{R}^k $ to $ \mathbb{R}^n $. A vector
$ \tau \in \mathbb{R}^k $ is called an $\varepsilon$-\emph{almost period} of $ g $ if
$$|g^{\tau}(x) -g(x)| < \varepsilon \quad \textrm{for all}\; x \in
\mathbb{R}^k.
$$

A set $ E \subset \mathbb{R}^k $ is called \emph{relatively dense} if there exists $L < \infty$ such that every
$k$-dimensional ball of radius $ L $ has a nonempty intersection with $ E$. It is obvious that we can replace
a ball with a cube in this definition.

A continuous mapping $ g $ from $ \mathbb{R}^k $ to $ \mathbb{R}^n $ is called
\emph{almost periodic} if for every $\varepsilon >0$  the set of $\; \varepsilon
$-almost periods of $ g $ is relatively dense in $ \mathbb{R}^k $.

Point out some properties of almost periodic mappings. First of all, the almost
periodicity of a mapping from $ \mathbb{R}^k $ to $ \mathbb{R}^n $ is equivalent to the almost
periodicity of its coordinates, therefore, it is sufficient to consider almost
periodic functions on $ \mathbb{R}^k $.

The following assertions are well known for almost periodic functions on the real axis
(see, for example, \cite{cord}, \cite{sib}). Their proofs don't change much in the
case of higher dimension (see, for example, \cite{ron tube}).

\begin{Pro}\label{Prop1}

a) An almost periodic function on $ \mathbb{R}^k $ is bounded and uniformly continuous,

b) if for an almost periodic function $ g $ the sequence $ (g^{\tau_n}) $ converges
uniformly on $ \mathbb{R}^k $ then its limit $ \widetilde{g} $ is an almost periodic
function as well. In addition, $g\not\equiv const $ if and only if $
\widetilde{g}\not\equiv const $.
\end{Pro}

\begin{Pro}\label{Prop2}
For a continuous mapping $ g(x) $ on $ \mathbb{R}^k $ the following conditions are equivalent:

a)  $ g(x) $ is almost periodic,

b) for each sequence $ (h_n)\subset\mathbb{R}^k $ there exists a subsequence $
(h_{n'}) $ such that $ |g(x+h_{l'})
-g(x+h_{m'})|\underset{{l'},{m'}\rightarrow\infty}{\longrightarrow} 0 $ uniformly on $
\mathbb{R}^k $,

c) there is a sequence of finite exponential sums
 $$
                S_n (x) = \sum_{i=1}^{N(n)} c_{i,n} e^{\langle\lambda_{i,n},x\rangle}
 $$
that converges to $ g(x) $ uniformly on $ \mathbb{R}^k $.
\end{Pro}

The assertion b)  is called the \emph{Bochner criterion} for almost periodic mappings.

{\bf 2. Almost periodic discrete multiple sets.}

As it was introduced in \cite{kolb}, we call a value set of a sequence
$(a_n)_{n=1}^\infty\subset\mathbb{R}^k$ a \emph{discrete multiple set} (we write
$\{a_n\}_{n\in\mathbb{N}}$) if this sequence doesn't possess any finite limit points.
In other words, a discrete multiple set is a discrete set where each point has a
finite multiplicity.

For any two discrete multiple sets $\{a_n\}_{n\in\mathbb{N}}$ and
$\{b_n\}_{n\in\mathbb{N}}$ we define a \emph{distance} between them be the formula
 $$
 dist(\{a_n\}_{n\in\mathbb{N}},\{b_n\}_{n\in\mathbb{N}}) =
\underset{\sigma}{\inf}\;\underset{n\in\mathbb{N} }{\sup}\,|\,a_n - b_{\sigma(n)}|,
 $$
where infimum is taken over all bijections $ \sigma
:\,\mathbb{N}\longrightarrow\mathbb{N} $. As was shown in \cite{kolb}, this function
satisfies all the axioms of metric except the finiteness.

\begin{Def}\label{Def1}
 A vector $\tau\in\mathbb{R}^k $ is called an $\varepsilon$-almost period of a discrete multiple set
 $\{a_n\}_{n\in\mathbb{N}}\subset\mathbb{R}^k$,
 if
  $$
  dist(\{a_n\}_{n\in\mathbb{N}},\{a_n+\tau\}_{n\in\mathbb{N}})<\varepsilon.
 $$
\end{Def}

\begin{Def}\label{Def2}
A discrete multiple set $\{a_n\}_{n\in\mathbb{N}}\subset\mathbb{R}^k$ is called almost
periodic, if for each $ \varepsilon>0$ the set of its $\varepsilon $-almost periods is
relatively dense in $\mathbb{R}^k $.
\end{Def}

Note that the sum and the difference of any two $\varepsilon $-almost periods are $2\varepsilon$-almost periods.

\begin{Th}\label{Th1}
The limit of almost periodic multiple sets is almost periodic as well.
\end{Th}
\bpr Let $(\{{a_n^{(p)}}\}_{n\in\mathbb{N}})$ be a sequence of almost periodic sets.
Let $\{b_n\}_{n\in\mathbb{N}} $ be such a set that
$dist\left(\{{a_n^{(p)}}\}_{n\in\mathbb{N}},\{b_n\}_{n\in\mathbb{N}}\right)\underset{p\rightarrow\infty}{\rightarrow}0$.
Let us take any $ \varepsilon>0 $. Let $ p_0$ be such a number that
$$
dist\left(\{ \{a_n^{(p_0)}\}_{n\in\mathbb{N}},\{{b_n}\}_{n\in\mathbb{N}}\right)< \frac{\varepsilon}{3}.
$$

\noindent Let $ E_\frac{\varepsilon}{3}$ be a relatively dense set of
$\frac{\varepsilon}{3} $-almost periods of $ \{{a_n^{(p_0)}}\}_{n\in\mathbb{N}} $. For
$ \tau\in E_\frac{\varepsilon}{3} $ we have
$$
dist\left(\{ {a_n^{(p_0)}}\}_{n\in\mathbb{N}},\{{a_n^{(p_0)}}+ \tau\}_{n\in\mathbb{N}}\right)< \frac{\varepsilon}{3}.
$$

\noindent Therefore,
 $$
dist\left(\{b_n \}_{n\in\N},\{ {b_n+ \tau} \}_{n\in\mathbb{N}}\right)
  \leq dist\left(\{ b_n \}_{n\in\N},\{{a_n^{(p_0)}}\}_{n\in\N}\right)
 $$
 $$
 + dist\left(\{{a_n^{(p_0)}}\}_{n\in\mathbb{N}},\{{a_n^{(p_0)}}+\tau\}_{n\in\mathbb{N}}\right)+
 dist\left(\{{a_n^{(p_0)}}+\tau
\}_{n\in\mathbb{N}},\{{b_n+ \tau}\}_{n\in\mathbb{N}}\right)< \varepsilon.
 $$
Hence, $\tau$ is an $ \varepsilon$-almost period of $\{b_n\}_{n\in\mathbb{N}} $, and
the set of almost periods is relatively dense. \bs

\noindent R e m a r k. As it was shown in \cite{kolb}, theorem 2, a space $ (X,dist)$
of all discrete multiple sets is complete \footnote{Precisely, a metric space of all
discrete multiple sets lying at a finite distance from some fixed discrete multiple
set is complete.}. Hence, almost periodic multiple sets form a complete closed
subspace in $(X,dist) $.

The following theorem is an analogue of the Bochner criterion for almost periodic
multiple sets.

\begin{Th}\label{Th2}
A discrete multiple set $\{a_n\}_{n\in\mathbb{N}}$ is almost periodic if and only if
for every sequence $(h_p)_{p=1}^\infty\subset \mathbb{R}^k$ there is a subsequence
$(h_p' )_{p=1}^\infty$ such that $(\{a_n + h_p'\}_{n\in\mathbb{N}})_{p=1}^\infty$ has
a limit.
\end{Th}

The proof of this theorem is based on the following lemma:

\begin{Le}\label{Le1} Let $D =  \{a_n\}_{n\in\mathbb{N}}$ be an almost periodic
multiple set. Then every sequence $(h_p)_{p=1}^\infty \subset \mathbb{R}^k$ has a
subsequence $(h_p ')_{p=1}^\infty$ such that for each $\varepsilon >0$ and arbitrary
numbers $ l, m > N=N(\varepsilon) $ one can find an $\varepsilon$-almost period $ \tau
$ of $D $ satisfying the inequality
$$|h_l'- h_m' - \tau|<\varepsilon.$$
\end{Le}
\bpr First put $\varepsilon = 1$. Let $E_\frac{1}{2}$ be the relatively dense set of
$\frac{1}{2}$-almost periods of $D $. Consider the sets $$A_p = \underset{\tau\in
E_\frac{1}{2}}{\bigcup}B(h_p+\tau,1/2),  p \in \mathbb{N} .$$

\noindent There exists $ L< \infty $ such that for all $ p=1,2,\ldots $ the ball $
B(-h_p,L) $ contains a $ \frac{1}{2}$-almost period $\,\tau_p $ of $ D $. Since $
\tau_p + h_p \in B(0,L)$, we have $ A_p\bigcap B( 0 ,L)\neq0 $ for all $ p=1,2,\ldots
$.

Let us show that there exists a point $ x \in B(0,L) $ belonging to an infinite
sequence of the sets $ A_p$. Cover the ball $ B(0,L) $ by a finite number of mutually
disjoint $k$-dimensional cubes with edges of length $
\displaystyle\frac{\varepsilon}{2k}$. The Dirichlet principle implies that there is a
cube containing an infinite sequence of points $ h_p+\tau $. The diagonal of such a
cube is less than $ \displaystyle \frac{1}{2\sqrt{k}} $, hence this cube is contained
in the balls $B(h_p+\tau,\frac{1}{2}) $, therefore, in an infinite sequence of the
sets $ A_p $. Hence there exists a subsequence $(A_p^{(1)})_{p=1}^\infty \subset
(A_p)_{p=1}^\infty $ such that
$$A_p^{(1)}  = \underset{\tau\in
E_\frac{1}{2}}{\bigcup}B(h_p^{(1)}+\tau,1/2), \; p\in \mathbb{N} $$
and
$$
\bigcap_{p=1}^\infty A_p^{(1)}  \neq 0.
$$

Take a number $  h\in A_i^{(1)} \bigcap A_j^{(1)} $ with some $i,j$. There exists $
\tau' \in E_\frac{1}{2} $ such that $ |h-(h_i^{(1)}+\tau ')|<\frac{1}{2} $ and $
\tau'' \in E_\frac{1}{2}$ such that $|h-(h_j^{(1)}+\tau '')|<\frac{1}{2}. $ Therefore
we have

$$
|h_i^{(1)}-h_j^{(1)}-(\tau''-\tau')| \leq |h-(h_j^{(1)}+\tau'')|+|h-(h_i^{(1)}+\tau')|< 1.
$$

Thus, for each $ i, j $ there exists $ 1 $-almost period $ \tau = \tau'' -\tau'$ of $D $ such that the inequality
$$
|h_i^{(1)}-h_j^{(1)}-\tau|< 1
$$
holds.

Now put $ \varepsilon =\frac{1}{2} $. Similarly, construct the sets
$$ A_p^{(2)} = \underset{\tau\in
E_\frac{1}{4}}{\bigcup}B(h_p^{(2)}+\tau,1/4), \; p \in \mathbb{N}$$ and
$(h_p^{(2)})_{p=1}^\infty\subset (h_p^{(1)})_{p=1}^\infty $ such that for all $ i, j $
the inequality
$$
|h_i^{(2)}-h_j^{(2)}-\tau|< \frac{1}{2}
$$
holds for some $ \frac{1}{2} $-almost period $\tau $ of $D$.

Repeating this construction for $ \varepsilon =\frac{1}{3}, \varepsilon =\frac{1}{4},
\ldots $ and choosing a diagonal subsequence $ (h_p ')_{p=1}^\infty $, we obtain the
assertion of the lemma. \bs

\vspace{10pt} \noindent{$\,$P r o o f $\,$ o f $\,$ t h e $\,$ t h e o r e m
\ref{Th2}.} Suppose that for every sequence $ (h_p)_{p=1}^\infty\subset\mathbb{R}^k $
there is a subsequence $(h_p')_{p=1}^\infty$ such that a sequence $(\{a_n + h_p'
\}_{n\in\mathbb{N}})_{p=1}^\infty$ has a limit. If $\{a_n\}_{n\in\mathbb{N}} $ is not
almost periodic, then for some $ \varepsilon_0 >0$ there exists a sequence of $
k$-dimensional balls $ (B_p)_{p=1}^\infty $ with infinitely increasing diameters $ l_p
$, such that no ball contains any $\varepsilon_0$-almost period of
$\{a_n\}_{n\in\mathbb{N}} $.

Let us take an arbitrary $ h_1 \in \mathbb{R}^k $ and a number $ \nu_1 $ such that $
\;l_{\nu_1}\,>\,1.$ For some $h_2 \in \mathbb{R}^k $ the difference $ h_2 - h_1 $
belongs to $ B_{\nu_1} $. Let $ \nu_2 $ be the first number such that $ \; l_{\nu_2}
\,>\,\max\{2,|h_2 - h_1|\}.$ Take $ h_3 \in \mathbb{R}^k $ such that the differences $
h_3 -h_1, h_3 - h_2 $ belong to $ B_{\nu_2} $ (it is possible, since the latter
condition is equivalent to $ | h_2 - h_1 |<l_{\nu_2}$). Generally, take a number
$\nu_n $ such that $ \; l_{\nu_n} \,>\,\max\{ n,|h_2 - h_1|,|h_3 - h_2|, \ldots, |h_n
- h_2|, |h_n - h_1|\}$ and $ h_{n+1} \in \mathbb{R}^k $ such that the differences
$h_{n+1} - h_1, h_{n+1} - h_2,\ldots, h_{n+1} - h_n $ belong to $ B_{\nu_n} $.

Take arbitrary $ p, m \,(p\,>\,m)$. By construction $ h_p- h_m \in\,B_{\nu_{p-1}} $,
and $B_{\nu_{p-1}}$ doesn't contain any $ \varepsilon_0$-almost period of
$\{a_n\}_{n\in\mathbb{N}} $. We have
$$dist(\{a_n+h_p\}_{n\in\mathbb{N}},\{a_n+h_m\}_{n\in\mathbb{N}}) = dist(\{a_n+(h_p -
h_m)\}_{n\in\mathbb{N}},\{a_n\}_{n\in\mathbb{N}}) \geq \varepsilon_0.$$ Thus, there
are no convergent subsequences in the sequence $(\{a_n +
h_{p}\}_{n\in\mathbb{N}})_{p=1}^\infty$.

Conversely, let a discrete multiple set $\{a_n\}_{n\in\mathbb{N}} $ be almost
periodic. Consider an arbitrary sequence $(h_p)_{p=1}^\infty \subset \mathbb{R}^k$. By
Lemma \ref{Le1} there exists a subsequence $ (h_p')_{p=1}^\infty$ such that for every
$\varepsilon >0$ and arbitrary $ l, m > N=N(\varepsilon) $ one can find an
$\varepsilon/2$-almost period $ \tau $ of $\{a_n\}_{n\in\mathbb{N}} $ with the
property
$$|h_l' -h_m' -\tau|<\frac{\varepsilon}{2}.$$

We get
\begin{gather*}
dist(\{a_n+h_l'\}_{n\in\mathbb{N}},\{a_n+h_m'\}_{n\in\mathbb{N}}) = \underset{bij\; \sigma
:\,\mathbb{N}\longrightarrow\mathbb{N}}{\inf}\;\underset{n\in\mathbb{N} }{\sup}\,|\,a_n +h_l' - a_{\sigma(n)}-
h_m'| \\
\leq\underset{bij\; \sigma
:\,\mathbb{N}\longrightarrow\mathbb{N}}{\inf}\underset{n\in\mathbb{N} }{\sup}\,|\,a_n
- a_{\sigma(n)}+\tau| +|h_l' - h_m' - \tau|\\<
dist(\{a_n\}_{n\in\mathbb{N}},\{a_{\sigma(n)}-\tau\}_{n\in\mathbb{N}})+\varepsilon/2<
\varepsilon.
\end{gather*}
Since $\varepsilon$ is arbitrary, we see that the sequence $(\{a_n +
h_{p}'\}_{n\in\mathbb{N}})_{p=1}^\infty$  is fundamental. Using the remark after
Theorem \ref{Th1}, we finish the proof. \bs

\begin{Def}\label{Definition S-property}
\emph{(\cite{kolb})} A discrete multiple set
$\{a_n\}_{n\in\mathbb{N}}\subset\mathbb{R}^k$ possesses $S$-property, if there exists
$ L <\infty $ such that for any $\tau\in\mathbb{R}^k $ there is a bijection $\sigma
:\,\mathbb{N}\longrightarrow\mathbb{N}$ with the property
$$
         \underset{n\in\mathbb{N}}{\sup}|(a_n+\tau)-a_{\sigma(n)}|  \leq L.
$$
\end{Def}

In the other words, a discrete multiple set
$\{a_n\}_{n\in\mathbb{N}}\subset\mathbb{R}^k$ possesses $S$-property, iff there exists
$ L <\infty $ such that the inequality $dist(\{a_n\}_{n\in\mathbb{N}},\{a_n+
\tau\}_{n\in\mathbb{N}}) \leq L $ holds for any $\tau\in\mathbb{R}^k $.

By \cite{kolb}, any almost periodic discrete multiple set possesses $S$-property.
Therefore Theorem 4 and Proposition 4 from \cite{kolb} imply the following statements:

\begin{Th}\label{Th3}
Let $D$ be an almost periodic multiple set. Then there exists $ M < \infty $ such that
\begin{equation}
       \label{1}
         \mathrm{card}\; (D\cap B(c,1))< M \quad \textrm{ for all } c\in \mathbb{R}^k.
\end{equation}
\end{Th}

\begin{Th}\label{Th4}
For any almost periodic multiple set $D$ there exists $C < \infty $ such that for any
convex bounded set $E\subset\mathbb{R}^k$ and $t\in\mathbb{R}^k$ the inequality
$$
  |\mathrm{card} \;(D\cap E) -\mathrm{card} \;(D\cap(E+t))|<C((\mathrm{diam}
  \;E)^{k-1}+1)
$$
is fulfilled.
\end{Th}

Following \cite{kolb}, the density of a discrete multiple set is the value
$$
\Delta = \underset{T\rightarrow\infty}{\lim}\displaystyle \frac{\mathrm{card} \;(D\cap
Q(0,T))}{ T^k }.
$$

From Theorem 5 of \cite{kolb} we get
\begin{Th}\label{Th5}
Any almost periodic multiple set possesses finite nonzero shift invariant density.
\end{Th}

\vspace{10pt}

{\bf 3. The connection between almost periodic measures and almost periodic discrete
multiple sets.}

\begin{Def}\label{Def}
\emph{(\cite{ron tube})} A locally finite complex-valued Radon measure $ \mu $ on
$\mathbb{R}^k$ is called almost periodic if for each compactly supported continuous
function $ \varphi $ on $\mathbb{R}^k $ the convolution
$$\varphi*\mu(z) = \int_{\mathbb{R}^k} \varphi (y-z)d\mu(y) $$
is almost periodic on $\mathbb{R}^k$.
\end{Def}

In other words, a locally finite complex-valued Radon measure $ \mu $ is almost
periodic if for each compactly supported continuous function $ \varphi $ on
$\mathbb{R}^k $ and any $ \varepsilon>0 $ there exists a relatively dense set
$E_{\varepsilon, \varphi}$ (the set of ($\varepsilon,\varphi$)-almost periods of the
measure $ \mu $) with the property
$$
|\varphi*\mu(z)-\varphi*\mu(z-\tau)|\leq\varepsilon\quad \forall \tau \in E_{\varepsilon, \varphi},\;
\forall z \in \mathbb{R}^k.
$$

Note that the sum and the difference of any two ($\varepsilon,\varphi$)-almost periods
are $(2\varepsilon,\varphi)$-almost periods.

We will say that measures $ \mu_n $ converge \emph{weakly uniformly} to some measure $ \mu $ if
$$
\varphi*\mu_n(x)\underset{n\rightarrow\infty}{\rightarrow}\varphi*\mu(x)
$$
uniformly on $ \mathbb{R}^k $ for every continuous compactly supported function $\varphi $.

Note some properties of almost periodic measures.

\begin{Th}\label{Th6}
If almost periodic measures converge weakly uniformly to some measure $ \mu $, then $ \mu $ is almost periodic as well.
\end{Th}

\noindent The proof follows immediately from the definition of an almost periodic
measure and Proposition \ref{Prop1} of the present article.

\begin{Th}\label{Th7}
Let $\mu $ be an almost periodic measure. Then there exists $ M < \infty $ such that the condition
$$
           |\mu|(B(c,1))\leq M \quad \textrm{ for all } c\in \mathbb{R}^k
$$
is fulfilled; here $ |\mu|$ is the variation of measure $ \mu$.
\end{Th}
\noindent This fact follows from Theorem 2.1. in \cite{ron tube} with $S=\{0\},\; p=0$.

The following theorem is an analogue of the Bochner criterion for almost periodic mappings.

\begin{Th}\label{Th8}
A Radon measure $ \mu $ is almost periodic if and only if for every sequence $
(h_n)\subset \mathbb{R}^k $ there is a subsequence $(h_{n'})$ such that measures
$\mu^{h_{n'}} $ converge weakly uniformly to some measure $\mu' $.
\end{Th}
\noindent This theorem is a corollary of Theorem 2.2. in \cite{ron tube} with $S=\{0\},\; p=0$.

\begin{Def}\label{Def5}
For any Radon measure $ \mu $ the value
$$
           \Delta = \underset{T\rightarrow\infty}{\lim}\displaystyle \frac{\mu( Q( 0,T))}{ T^k }
$$
is called the density of $ \mu $.
\end{Def}

\begin{Th}\label{Th9}
Any almost periodic Radon measure possesses finite shift invariant density, i.e.,
 $$
 \Delta =\underset{T\rightarrow\infty}{\lim}\displaystyle \frac{\mu( Q( \alpha,T))}{ T^k }
 $$
uniformly over $\alpha\in \mathbb{R}^k$.
\end{Th}
\noindent (See \cite{ron tube}, Theorem 2.7. with $S=\{0\}. $)

We can associate with a discrete multiple set $D = \{ a_n\}_{n\in\mathbb{N}}$ the
measure
$$
\mu_D =   \sum_{ a_n \in D} \delta(x-a_n). 
$$
Note that for every continuous function $\varphi $
$$
\varphi*\mu_D(x) =  \sum_{n\in\mathbb{N}} \varphi(x-a_n),\; x\in\mathbb{R}^k.
$$

\begin{Th}\label{Th10}
Let  $(D_p) = (\{ a_n^{(p)} \}_{n\in\mathbb{N}})$ be a sequence of discrete multiple
sets. Let $D=\{a_n\}_{n\in\mathbb{N}}$ be a discrete multiple set satisfying
({\ref{1}}). Then the following conditions are equivalent:

a) discrete multiple sets $ D_p $ converge to $ D $,

b) measures $ \mu_{D_p} $ converge weakly uniformly to $\mu_D$.
\end{Th}

We need the following lemma:

\begin{Le}\label{Le2}
Let  $(D_p) = (\{ a_n^{(p)} \}_{n\in\mathbb{N}})$ be a sequence of discrete multiple
sets, which converges to a discrete multiple set $D=\{a_n\}_{n\in\mathbb{N}}$
satisfying ({\ref{1}}). Then for sufficiently large $p$ any discrete multiple set
$D_p$ satisfies ({\ref{1}}) (with a constant $4^kM$ instead of $M$).
\end{Le}
\bpr For sufficiently large $p$ we have
$\underset{\sigma}{\inf}\;\underset{n\in\mathbb{N} }{\sup}\,|\,a_n -
a_{\sigma(n)}^{(p)}| <1 $. Hence, $\mathrm{card} \; (D_p \cap B(c,1))\leq
\mathrm{card} \; (D \cap B(c,2)) \leq 4^kM$ for all $c \in \mathbb{R}^k$. \bs

\vspace{10pt} \noindent{ $\,$P r o o f $\,$ o f $\,$ t h e o r e m \ref{Th10}.} Let
the measures $ \mu_{D_p} $ converge weakly uniformly to $\mu_D$. We will show that $
dist\left(\{
{a_n^{(p)}}\}_{n\in\mathbb{N}},\{a_n\}_{n\in\mathbb{N}}\right)\underset{p\rightarrow\infty}{\rightarrow}0
$, i.e. for any $ \varepsilon>0 $ for each sufficiently large $ p $ there exists a
bijection $\sigma:\,\mathbb{N}\longrightarrow\mathbb{N}$ such that
\begin{equation}
       \label{2}
       \underset{n\in\mathbb{N} }{\sup}\,|\,a_n^{(p)} - a_{\sigma(n)}|< \varepsilon.
\end{equation}

First note that for any $ \varepsilon>0 $ there is $\eta>0$ such that the diameter of
every connected component of the union $ \bigcup_n B(a_n,\eta)$ is less than $
\varepsilon $. Let us check it. Put $\varepsilon\in(0,\,1)$ and put $\eta<
\varepsilon/(2M+1)$, where $M$ satisfies ({\ref{1}}). Let $A$ be an arbitrary
connected component of the union $ \bigcup_n B(a_n,\eta)$ and $c$ be an arbitrary
point of $A$. By Theorem \ref{Th3}, the ball $B(c,1)$ contains at most $M$ points of
$D$. Choose from them the points $a_{n_1}, \ldots,  a_{n_N} \;(N\leq M)$ belonging to
$ A$ (we suppose that $B(a_{n_i},\eta)\cap B(a_{n_{i+1}},\eta)\neq 0, \;
i=\overline{1,N-1}$). If $A$ isn't contained in $B(c,1)$, then there exists a point
$a'\in A\setminus B(c,1)$ with the property $|a' - a_{n_i}|< 2\eta$ for some
$n_i,\,i=\overline{1,N}$. Since $ c \in B(a_{n_p},\eta)$ for some
$n_p,\,p=\overline{1,N}$, we have
$$|a' - c|\leq |a' -a_{n_i}| + (N-1)\underset{j=\overline{1,N-1}}{\max}|a_{n_j}-a_{n_{j+1}}|+|c-a_{n_p}|
$$
 $$
 < 2\eta +(N-1)2\eta +\eta  = \eta (2N+1)<1,
 $$
that is impossible. Thus $A\subset B(c,1)$ and
 $$
\textrm{diam} A  < 2\eta M< \varepsilon.
 $$

Take a nonnegative continuous function $ \varphi $ with the support in $B(0,\eta/2) $
such that $ 0\leq\varphi(x)\leq \varphi(0) = 1 $. Put $ \nu  = \int\varphi dm$. We may
assume that $\nu < 1$. Since measures $ \mu_{D_p} $ converge weakly uniformly to
$\mu_D$, for sufficiently large $p$ we have
\begin{equation}
\label{3}|\varphi* \mu_{D_p}(x)-\varphi*\mu_D(x)|< \nu/2\quad \forall\,
x\in\mathbb{R}^k.
\end{equation}

\noindent The distance between an arbitrary pair of terms $ a_n \in A $ and $ a_m
\notin A $ is at least $ 2 \eta $, therefore,
$$
\int_A\left(\int_{\mathbb{R}^k} \varphi(x - y) d\mu_D (y) \right)dm(x) =\int_A\left(\int_A \varphi(x - y) d\mu_D
(y) \right)dm(x)
$$
\begin{equation}
\label{4}=\int_A\left(\int_{\mathbb{R}^k} \varphi(x - y)dm(x) \right)d\mu_D (y) = \nu
\:\mathrm{card} \; (D\cap A).
\end{equation}

Furthermore, let $ a_n^{(p)} $ be a term of $ D_p $ such that $ a_n^{(p)}  \in A$. We
have $ \varphi * \mu_{D_p} (a_n^{(p)} ) \geq \varphi(0) = 1. $ In view of ({\ref{3}}),
we get
$$
\varphi * \mu_{D} (a_n^{(p)} ) \geq 1-\nu/2 >1/2.
$$

\noindent Hence, there exists a term $ a_{n'} $ of $ D $ such that $ | a_{n'} -
a_n^{(p)} |<\eta/2$. Then the distance between $ a_{n'} $ and $ A $ is less than $
\eta/2 $. Therefore, $ a_{n'} \in A  $. Similarly, for any $ a_j^{(p)} \notin A $
there exists a term $ a_{j'} $ of $ D $ such that $ | a_{j'} - a_j^{(p)} |<\eta/2$.
Consequently, the distance between $ a_{j'} $ and $\mathbb{R}^k \backslash A $ is less
than $ \eta/2 $, so  $ a_{j'} \notin A $. Therefore,
$$
\int_A\left(\int_{\mathbb{R}^k} \varphi(x - y) d\mu_{D_p}(y) \right)dm(x) =\int_A\left(\int_A \varphi(x - y)
d\mu_{D_p}(y) \right)dm(x)
$$
\begin{equation}
\label{5}=\int_A\left(\int_{\mathbb{R}^k} \varphi(x - y)dm(x) \right)d\mu_{D_p}(y) =
\nu \:\mathrm{card} \; (D_p \cap A).
\end{equation}

\noindent Since
$$
\left|\int_A(\varphi* \mu_{D_p}(x)-\varphi*\mu_D(x)) dm(x)\right|\leq m(A) \:\nu/2 <\nu/2,
$$
the values from equalities ({\ref{4}}) and ({\ref{5}}) coincide. Therefore, $
\mathrm{card}\; (D \cap A)= \mathrm{card} \; (D_p \cap A)$. The same is valid for all
connected components of union $ \bigcup_n B(a_n,\eta)$. On the other hand, we have
just proved that each term $ a_n^{(p)}  $ of $ D_p $ belongs to $B(a_{n'},\eta) $ for
some term $ a_{n'} $ of $ D $. Hence, it belongs to a connected component of $
\bigcup_n B(a_n,\eta)$. Consequently, there is a bijection $ \sigma $ such that
({\ref{2}}) is fulfilled.

Conversely, suppose that discrete multiple sets $D_p$ converge to a discrete multiple
set $ D $. Let $ \varphi $ be a nonnegative continuous function with compact support $
G $ in $ \mathbb{R}^k $. We will show that the convolutions $\varphi* \mu_{D_p}(x)$
converge to  $\varphi*\mu_D(x)$ as $p\rightarrow\infty $ uniformly over $ x \in
\mathbb{R}^k $.

Take $ z \in \mathbb{R}^k$ and $\varepsilon >0$. There is $ \delta >0$ such that the
inequality $ |x_1-x_2|<\delta $ implies $|\varphi(x_1) - \varphi(x_2)|< \varepsilon $.
By assumption, for each sufficiently large $ p $ there exists a bijection $\sigma$
such that $ |a_n - a_{\sigma(n)}^{(p)}|<\delta$ for all $ n\in \mathbb{N} $. By Lemma
\ref{Le2}, the set $ G $ contains at most $ M' = M'(G) < \infty $ terms of $ D_p $.
Hence,
$$
|\varphi* \mu_{D_p}(x)-\varphi*\mu_D(x)|\leq M' \underset{|x_1-x_2|<\delta}{\sup}|\varphi(x_1) - \varphi(x_2)| <M'
\varepsilon.
$$
\bs

Now we can prove the main result of this section.

\begin{Th}\label{Th11}
A discrete multiple set $ D $ is almost periodic if and only if the corresponding
measure $ \mu_D $ is almost periodic.
\end{Th}
\bpr By theorem \ref{Th2} a discrete multiple set $ D $ is almost periodic if and only
if for every sequence $(h_p)_{p=1}^\infty\subset \mathbb{R}^k$ there is a subsequence
$(h_p')_{p=1}^\infty$ such that the sequence of almost periodic multiple sets $(\{a_n
+ h_p'\}_{n\in\mathbb{N}})_{p=1}^\infty$ has a limit. By theorem \ref{Th10}, it is
true if and only if for every sequence $(h_p)_{p=1}^\infty\subset \mathbb{R}^k$ there
is a subsequence $(h_p')_{p=1}^\infty$ such that the sequence of almost periodic
measures $(\mu_{D + h_{p}'})_{p=1}^\infty$ converges weakly uniformly to some measure.
By theorem \ref{Th8} it is true if and only if the measure $\mu$ is almost periodic.
\bs

\end{document}